\documentclass[11pt,english]{article}

\usepackage[english]{babel}
\usepackage{amsthm}
\usepackage{amsmath}
\usepackage{amssymb}
\usepackage{graphicx}
\usepackage{dsfont}
\usepackage[bf]{caption}
\usepackage{hyperref} 
\usepackage{parskip}

\usepackage{anysize}
\marginsize{2.5cm}{2.5cm}{2cm}{1.9cm}

\newtheorem{thm}{Theorem}[section]
\newtheorem{cor}[thm]{Corollary}
\newtheorem{lem}[thm]{Lemma}
\newtheorem{prop}[thm]{Proposition}

\theoremstyle{definition}

\newtheoremstyle{rmk}
  {12pt}                   
  {12pt}                   
  {}                       
  {}                       
  {\normalfont\bfseries}   
  {.}                      
  { }               
  {}

\theoremstyle{rmk}

\newtheorem{rmk}[thm]{Remark}

\renewcommand{\qedsymbol}{$\blacksquare$}

\newcommand{\p} {\textnormal{\textsf{P}}}
\newcommand{\E} { \textnormal{\textsf{E}}}
\newcommand{\N} { \mathbb{N} }
\newcommand{\Z} { \mathbb{Z} }
\newcommand{\R} { \mathbb{R} }

\newcommand{\var}{\mathop{\mathsf{Var}}}

\newenvironment{pfof}[1]{\noindent{\underline{\textit{Proof of #1.}}}\hspace{0.5em}}
	{\hfill\qed\vspace{1ex}}

\allowdisplaybreaks

\begin{document}

\title{On the Strong Recurrence of Recurrent RWRE}
 \author{Michael Kochler}
\date{\today}
\maketitle

\begin{abstract}
\noindent We consider a recurrent RWRE $(X_n)_{n \in \N_0}$ on $\Z$ and investigate the quenched return probabilities of the RWRE to the origin for which we state results on their decay in terms of summability. Additionally, we give some examples for recurrent RWRE with a multidimensional state space which give reason for the part ``strong recurrence'' in the title of this paper when we compare the behaviour of RWRE with the behaviour of the symmetric random walk on $\Z^d$.\\


\textbf{AMS 2000 Subject Classiﬁcation:} Primary 60K37, 60J10
\end{abstract}

\section{Introduction}
In \cite{CP}, Comets and Popov also consider the return probabilities of the one-dimensional recurrent RWRE on $\Z$. In contrast to our setting, they consider the corresponding jump process in continuous time $(\xi^x_t)_{t \ge 0}$ started at $x \in \Z$ and with jump rates $(\omega_x^+,\omega_x^-)_{x \in \Z}$ to the right and left neighbouring sites. One advantage of this process in continuous time is that it is not periodic as the RWRE in discrete time. 
As one result, they show the following (under two conditions on the environment $(\omega_x^+,\omega_x^-)_{x \in \Z}$):

\textbf{Theorem} (cf.\ Corollary 2.1 and Theorem 2.2 in \cite{CP}) \textit{We have}
\[
\frac{\log P_{\omega}(\xi_t^0 = 0)}{\log t} \xrightarrow{t \to \infty} - \widehat{a}_e
\]
\textit{in law, where $\widehat{a}_e$ has the density}
\[
p(z)=\begin{cases}
2 - z - (z+2) \cdot e^{-2z} & \text{if } z \in (0,1) \\
\big([e^2-1] \cdot z - 2 \big) \cdot e^{-2z} & \text{if } z  \ge 1.
\end{cases}
\]
Since we can embed the recurrent RWRE $(X_n)_{n \in \N_0}$ in discrete-time into the corresponding jump process in continuous time, we can expect the return probabilities to behave similarly as in the continuous setting. In particular, for $\p$-a.e.\ environment $\omega$, we expect 
\begin{align*} \label{preI.1}
P_{\omega}(X_{2n} = 0) =: n^{-a(\omega,n)} 
\end{align*}
with
\begin{align*}
& \liminf_{n \to \infty} a(\omega,n) = 0,\ \limsup_{n \to \infty} a(\omega,n) = \infty.
\end{align*}

In order to answer the question on recurrence and transience of particular examples of RWRE with different state spaces, one needs to know more about the decay of $(P_{\omega}(X_{2n}=0))_{n \ge 0}$ for fixed environment $\omega$. For the examples included in Section \ref{I-sec1.6}, the following two statements will be helpful (cf.\ Theorem \ref{I-Rthm1} and \ref{I-Rthm2}): For $\p$-a.e.\ environment $\omega$, we have
\begin{alignat*}{4}
& \sum_{n \in \N} P_{\omega}(X_{2n}=0) \cdot n^{-\alpha}= \infty &\qquad & \text{for } 0 \le \alpha < 1 \text{ and} & \\ 
& \sum_{n \in \N} \Big(P_{\omega}(X_{2n}=0)\Big)^{\alpha} = \infty &\qquad  & \text{for } \alpha > 0.&
\end{alignat*}
The structure of this paper is the following: In Section \ref{I-sec1.2}, we introduce the model of a RWRE on $\Z$ together with the notation which we use in this paper. Then we collect some useful equalities and inequalities in the context of RWRE in Section \ref{I-sec1.3} before we state our main result in Section \ref{I-sec1.4}. Section \ref{I-sec1.5} contains the proofs of our main results. The main tool for our proofs is a careful analysis of the corresponding potential of the RWRE  
(cf.\ \eqref{I-RWRE.a}). To this end, we introduce favourable ``valleys'' (cf.\ Figure \ref{I-figure1} on page \pageref{I-figure1}), which help us to derive lower bounds for the quenched return probabilities of the RWRE to the origin. In the last Section \ref{I-sec1.6}, we give some examples for recurrent RWRE with a multidimensional state space. In particular Corollary \ref{I-cor5.3a} and Corollary \ref{I-cor5.3} give reason for the part ``strong recurrence'' in the title of this paper when we compare the behaviour of RWRE with the behaviour of a symmetric random walk on $\Z^d$.

\section{Model and Notation} \label{I-sec1.2}
Let us first introduce the notation for a one-dimensional random walk in random environment (RWRE):\\
At first, let $\omega = (\omega_x)_{x \in \Z}$ be a sequence of i.i.d.\ random variables taking values in $(0,1)$ with respect to some probability measure $\p$. For $i \in \Z$ we define
\[
\rho_i = \rho_i(\omega):= \frac{1-\omega_i}{\omega_i}.
\]
In the following, we will assume that
\begin{align}
& \E[ \log \rho_0]= 0, \label{I-ass1} \vphantom{\int}\\
& \p(\varepsilon \le \omega_0 \le 1 -\varepsilon) = 1 \text{ for some } \varepsilon \in \left(0, \tfrac12 \right), \label{I-ass2} \vphantom{\int}\\
& \var (\log \rho_0) > 0. \label{I-ass3} \vphantom{\int}
\end{align}
Here, \eqref{I-ass1} ensures that the RWRE is recurrent. The second assumption is a common technical condition in the context of RWRE. Further, the third assumption excludes the case of a symmetric random walk on $\Z$.

For each environment $\omega$, we can introduce the random walk $(X_n)_{n \in \N_0}$ whose transition probabilities are determined by $(\omega_x)_{x \in \Z}$. More precisely for every $x \in \Z$, $(X_n)_{n \in \N_0}$ is a Markov chain with respect to $P_{\omega}^{x}$ determined by
\begin{align}
& P_{\omega}^{x}(X_0=x) = 1 \nonumber \vphantom{\int},\\
& P_{\omega}^{x}(X_{n+1} = y+1|X_n=y) = \omega_y = 1 - P_{\omega}^{x}(X_{n+1} = y-1|X_n=y) \quad \forall y \in \Z. \label{I-RWRE} \vphantom{\int}
\end{align}
As usual, we use $P_{\omega}^{o}$ instead of $P_{\omega}^{0}$ and will even drop the superscript $o$ where no confusion is to be expected. We can now define the potential $V$ as
\begin{align} \label{I-RWRE.a}
V(x):= \begin{cases} \sum\limits_{i=1}^x \log \rho_i & \text{for }x=1,2,\ldots \\ 0 & \text{for }x=0 \\ \sum\limits_{i=x+1}^{0} \log (\rho_i)^{-1} & \text{for }x=-1,-2,\ldots\ . \end{cases}
\end{align}
Note that $V(x)$ is a sum of i.i.d.\ random variables which are centred and which are bounded by $C:=\log (1-\varepsilon) - \log \varepsilon > 0$ due to the assumptions \eqref{I-ass1} and \eqref{I-ass2}. One of the crucial facts for the RWRE is that, for fixed $\omega$, the random walk is a reversible Markov chain and can therefore be described as an electrical network. The conductances are given by 
\[
C_{(x,x+1)}(\omega) = e^{-V(x)} = \begin{cases}  \prod\limits_{i=1}^{x} (\rho_i)^{-1} & \text{for }x=1,2,\ldots \\ 1 & \text{for }x=0 \\ \prod\limits_{i=x+1}^{0} \rho_i & \text{for }x=-1,-2,\ldots\end{cases}
\]
and the reversible measure which is unique up to multiplication by a constant is given by
\begin{align} \label{I-eq1.1}
\mu_{\omega}(x)= e^{-V(x)} + e^{-V(x-1)} = \begin{cases}  \prod\limits_{i=1}^{x-1} \frac{\omega_i}{1-\omega_i} \cdot \frac{1}{1-\omega_x} &  \text{for }x=1,2,\ldots \\ \frac{1}{\omega_0} & \text{for } x=0 \\    \prod\limits_{i=x+1}^{0} \frac{1-\omega_i}{\omega_i} \cdot \frac{1}{\omega_x} & \text{for }x=-1,-2,\ldots\ . \end{cases}
\end{align}
As a consequence of the reversibility, we conclude that we have
\begin{equation}
\label{I-eq1}
\mu_{\omega}(x) \cdot P^x_{\omega}(X_n = y) = \mu_{\omega}(y) \cdot P^y_{\omega}(X_n = x)
\end{equation}
for all $n \in \N_0$ and $x,y \in \Z$.
\section{Preliminaries}  \label{I-sec1.3}
In the following, we collect some useful properties of the RWRE. For the random time of the first arrival in $x$
\begin{equation}\label{I-def1}
\tau(x):= \inf \{n \ge 0:\ X_n=x\},
\end{equation}
the interpretation of the RWRE as an electrical network helps us to compute the following probability for $x < y < z$ (for a proof see for example formula (2.1.4) in \cite{Zei}):
\begin{equation}
\label{I-prel1}
P^y_{\omega}(\tau(z) < \tau(x)) = \frac{\sum\limits_{j=x}^{y-1} e^{V(j)}}{\sum\limits_{j=x}^{z-1} e^{V(j)}}
\end{equation}
Further (cf.\ (2.4) and (2.5) in \cite{SZ} and Lemma 7 in \cite{Gol}), we have for $k \in \N$ and $y < z$
\begin{align}
\label{I-prel2}
& P^y_{\omega}(\tau(z) < k) \le k \cdot \exp \left( - \max_{y \le i < z} \big[V(z-1)-V(i) \big] \right) 
\intertext{and similarly for $x < y$}
\label{I-prel3}
& P^y_{\omega}(\tau(x) < k) \le k \cdot \exp \left( - \max_{x < i \le y} \big[V(x+1)-V(i) \big] \right).
\end{align}
To get bounds for large values of $\tau(\cdot)$, we can use that for $x < y < z$ we have (cf.\ Lemma 2.1 in \cite{SZ})
\begin{equation}
\label{I-prel4}
E_{\omega}^y[\tau(z) \cdot \mathds{1}_{\{\tau(z) < \tau(x)\}}] \le (z-x)^2 \cdot \exp \left( \max_{x\le i \le j \le z} \big(V(j) - V(i)\big) \right).
\end{equation} 
Further, the Koml{\'o}s-Major-Tusn{\'a}dy strong approximation theorem (cf.\ Theorem 1 in \cite{KMT}, see also formula (2) in \cite{CP}) will help us to compare the shape of the potential with the paths of a two-sided Brownian motion:\\

\begin{thm} \label{I-Komlos}
In a possibly enlarged probability space, there exists a version of our environment process $\omega$ and a two-sided Brownian motion $(B(t))_{t \in \R}$ with diffusion constant $\sigma:=~(\var(\log \rho_0))^{\frac12}$ (i.e.\ $Var(B(t))=\sigma^2 |t|$) such that for some $K>0$ we have
\begin{equation}
\label{I-approx}
\p \left( \limsup_{x \to \pm \infty} \frac{|V(x)-B(x)|}{\log |x|} \le K \right) =1.
\end{equation} 
\end{thm}
\section{Results} \label{I-sec1.4}
Let us consider a RWRE $(X_n)_{n \in \N_0}$ on $\Z$ where the law of the environment $\omega=(\omega_x)_{x \in \Z}$ fulfils the assumptions $\eqref{I-ass1}, \eqref{I-ass2}$, and $\eqref{I-ass3}$. Then, the following two theorems hold:\\
\begin{thm} \label{I-Rthm1}
For $0 \le \alpha < 1$, we have
\begin{equation}
\label{I-thm1}
\sum_{n \in \N} P_{\omega}(X_{2n}=0) \cdot n^{-\alpha}= \infty
\end{equation}
for $\p$-a.e.\ environment $\omega$.
\end{thm}\vspace{12pt}

\begin{thm} \label{I-Rthm2}
For all $\alpha > 0$, we have
\begin{equation}
\label{I-thm2}
\sum_{n \in \N} \Big(P_{\omega}(X_{2n}=0)\Big)^{\alpha} = \infty
\end{equation}
for $\p$-a.e.\ environment $\omega$.
\end{thm}

\noindent For the last theorem we consider a combination of $d$ environments:\\

\begin{thm} \label{I-Rthm3}
For $d \in \N$, consider $d$ i.i.d.\ environments $\omega^{(1)}, \omega^{(2)}, \ldots, \omega^{(d)}$ which all fulfil the assumptions $\eqref{I-ass1}, \eqref{I-ass2}$, and $\eqref{I-ass3}$. Then, we have
\begin{equation}
\label{I-eqcor2}
\sum_{n \in \N} \prod_{k=1}^{d} P_{\omega^{(k)}}(X_{2n}=0) = \infty
\end{equation}
for $\p^{\otimes d}$-a.e.\ environment $(\omega^{(1)}, \omega^{(2)}, \ldots, \omega^{(d)})$.
\end{thm}

\begin{rmk}
A proof for Theorem \ref{I-Rthm3} can also be found in \cite{Zei} after Lemma A.2. The proof there uses the Nash-Williams inequality in the context of electrical networks.
\end{rmk}

\section{Proofs} \label{I-sec1.5}
Let us first introduce the sets $\Gamma^{+}(L,\delta)$ and $\Gamma^{-}(L,\delta)$ of environments for $L \in \N$ and \mbox{$0 < \delta < 1$} defined by
\begin{align*}
&\Gamma^{+}(L,\delta):=\{R^{+}_1(L) \le \delta L,\, R^{+}_2(L) \le \delta L,\, T^{+}(L) \le L^2\} \vphantom{\inf_{0 \le k \le T^{+}(L)}},\\
&\Gamma^{-}(L,\delta):=\{R^{-}_1(L) \le \delta L,\, R^{-}_2(L) \le \delta L,\, -T^{-}(L) \le L^2\} \vphantom{\inf_{0 \le k \le T^{+}(L)}},
\end{align*}
where
\begin{align*}
& T^{+}(L):= \inf \{n \ge 0:\ V(n) - \min_{0 \le  k \le n } V(k) \ge L\} \vphantom{\inf_{0 \le k \le T^{+}(L)}},\\
& T^{-}(L):= \sup \{n \le 0:\ V(n) - \min_{n \le k \le 0} V(k) \ge L\} \vphantom{\inf_{0 \le k \le T^{+}(L)}},\\
& R^{+}_1(L):= - \min_{0 \le k \le T^{+}(L)} V(k), \\
& R^{-}_1(L):= - \min_{T^{-}(L)\le k \le 0} V(k), \\
& T^{+}_b(L):= \inf \{n \ge 0:\ V(n) = - R^{+}_1(L)\}  \vphantom{\inf_{0 \le k \le T^{+}(L)}},\\
& T^{-}_b(L):= \sup \{n \le 0:\ V(n) = - R^{-}_1(L)\} \vphantom{\inf_{0 \le k \le T^{+}(L)}},\\
& R^{+}_2(L):= \max_{0 \le k \le T^{+}_b(L)} V(k), \\
& R^{-}_2(L):= \max_{T^{-}_b(L) \le k \le 0} V(k). \vphantom{\inf_{0 \le k \le T^{+}(L)}}
\end{align*}
Here, the $+$-sign and the $-$-sign indicate whether we deal with properties of the valley on the positive or negative half-line, respectively. Note that the definition of the sets $\Gamma^{+}(L,\delta)$ and $\Gamma^{-}(L,\delta)$ is compatible with the scaling of a Brownian motion in space and time (cf.\ \eqref{I-eq16a}).

\begin{figure}[ht] 
\vspace{45pt} \hspace{2cm}
\includegraphics [viewport=105 450 380 675, scale=1.05]{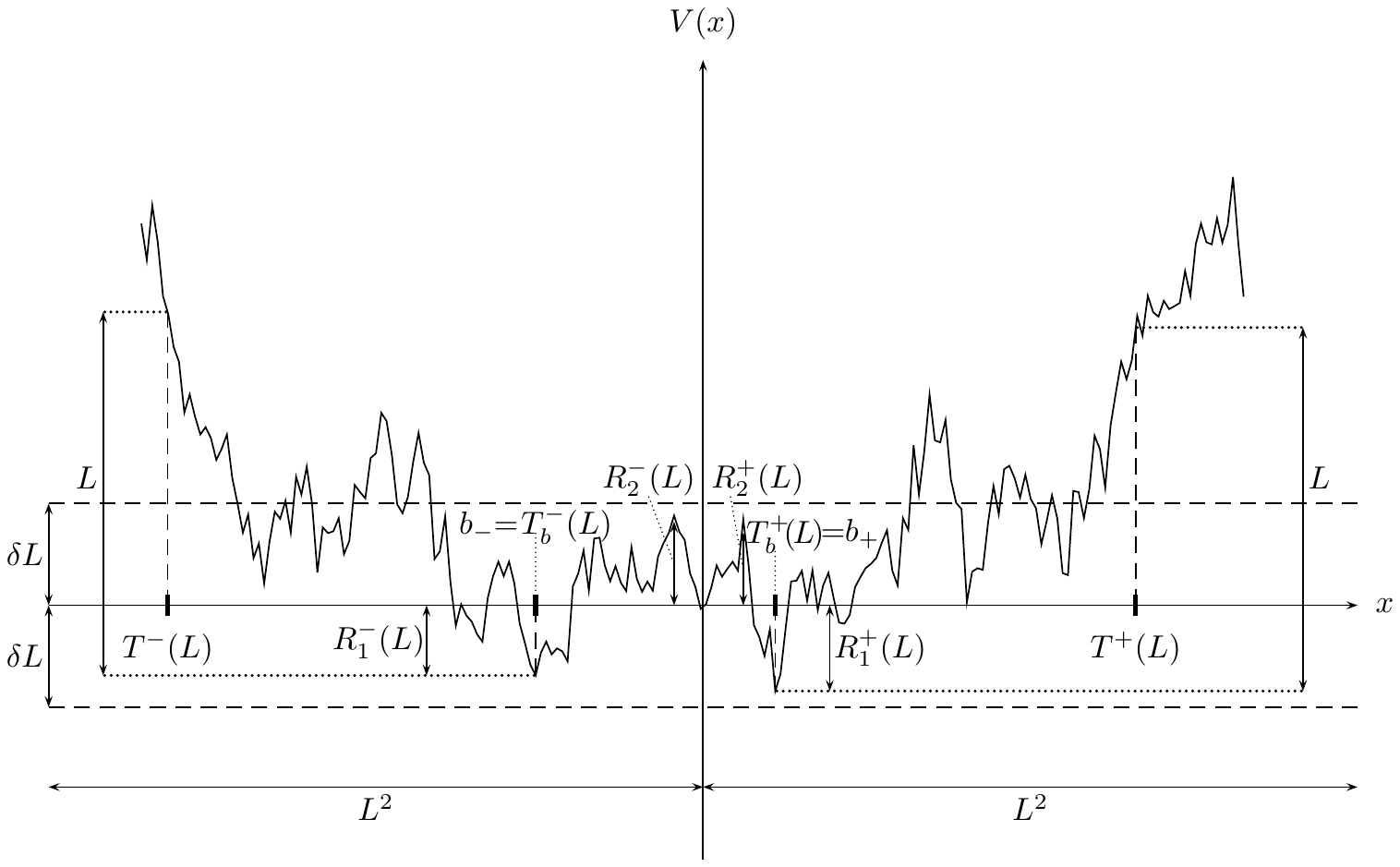} 
  \caption{Shape of a valley of an environment in $\Gamma(L,\delta):=\Gamma^{+}(L,\delta) \cap \Gamma^{-}(L,\delta)$} \label{I-figure1}
\end{figure}

\begin{rmk}
We have constructed the valleys in such a way that the return probability of the random walk to the origin is high (or bounded from below as we will see) for even time points as long as the random walk has not left the valley. For $\omega \in \Gamma^{+}(L,\delta) \cap \Gamma^{-}(L,\delta)$, we have the following behaviour for the random walk $(X_n)_{n \in \N_0}$ in the environment $\omega$:

\begin{enumerate}
\item Since we have $V(T^{-}(L)) - V(T^{-}_b(L)) \ge L$ and $V(T^{+}(L)) - V(T^{+}_b(L)) \ge L$, the random walk $(X_n)_{n \in \N_0}$ stays within $\{T^{-}(L), T^{-}(L) + 1, \ldots, T^{+}(L)\}$ for (approximately) at least $\exp(L)$ steps (cf.\ \eqref{I-eq4a}).

\item Within the area $\{T^{-}(L), T^{-}(L) + 1, \ldots, T^{+}(L)\}$, the random walk prefers to stay at positions $x$ with a small potential $V(x)$, i.e.\ at positions close to the bottom points at $T^{-}_b(L)$ and $T^{+}_b(L)$.

\item The return probability for the random walk from the positions $T^{-}_b(L)$ and $T^{+}_b(L)$ to the origin is mainly influenced by the potential differences $R_2^{-}(L)+R_1^{-}(L) \le 2 \delta L$ and $R_2^{+}(L)+R_1^{+}(L) \le 2 \delta L$ respectively, i.e.\ by the height of the climb the random walk has to trespass from the bottom points back to the origin (cf.\ \eqref{I-eq2}). 
\end{enumerate}

\end{rmk}

\begin{prop} \label{I-prop1}
For $\omega \in \Gamma(L,\delta):=\Gamma^{+}(L,\delta) \cap \Gamma^{-}(L,\delta)$ with $0 < \delta < \tfrac15$, we have
\begin{equation}
\label{I-lem1}
P_{\omega} (X_{2n}=0) \ge C \cdot \exp(-3\delta L)
\end{equation}
for \[ \exp (3 \delta L ) \le n \le \exp\big((1-2\delta)L\big),\]
where the constant $C=C(\delta)$ does not depend on $L$.
\end{prop}
\begin{pfof}{Proposition \ref{I-prop1}}
The construction of ``valleys'' has been useful for the proofs of many theorems in the context of RWRE. Our construction uses some ideas from \cite{CP}, where it is shown that the transition probabilities of a RWRE in continuous time converge in distribution. Since we deal with a RWRE in discrete time and we want to have lower estimates for the return probabilities for a fixed environment in Proposition \ref{I-prop1}, we will have to adapt the construction to our setting:\\[10pt]
The return probability to the origin for the time points of interest is mainly influenced by the shape of the ``valley'' of the environment $\omega$ between $T^-(L)$ and $T^+(L)$. For the positions of the two deepest bottom points of this valley on the positive and negative side, we write 
\[b_+:=T^+_b(L) \quad \text{and} \quad b_-:=T^-_b(L)\]
and we assume for the following proof that we have (cf.\ \eqref{I-def1} for the definition of $\tau(\cdot)$)
\begin{equation}
\label{I-ass4}
P^o_{\omega} \big( \tau(b_+) < \tau(b_-) \big) \ge \frac12.
\end{equation}
(Due to the symmetry of the RWRE, the proof also works in the opposite case if we switch the roles of $b_+$ and $b_-$). We have
\begin{align}
& P^o_{\omega} (X_{2n}=0) \ge P^o_{\omega} \left(X_{2n}=0,\ \tau(b_+) \le \tfrac{2n}{3},\ \tau(b_+) < \tau(b_-) \right) \nonumber \vphantom{\frac{\mu_{\omega}(0)}{\mu_{\omega}(b_+)}}\\
\ge\ & P^o_{\omega} \left(\tau(b_+) \le \tfrac{2n}{3},\ \tau(b_+) < \tau(b_-) \right) \cdot \widehat{\inf}_{\ell \in \big\{\left\lceil \tfrac{4n}{3}\right\rceil,\ldots,2n\big\}} P^{b_+}_{\omega} (X_{\ell} = 0) \nonumber \vphantom{\frac{\mu_{\omega}(0)}{\mu_{\omega}(b_+)}}\\
=\ & P^o_{\omega} \left(\tau(b_+) \le \tfrac{2n}{3},\ \tau(b_+) < \tau(b_-) \right) \cdot \frac{\mu_{\omega}(0)}{\mu_{\omega}(b_+)} \cdot \widehat{\inf}_{\ell \in \big\{\left\lceil \tfrac{4n}{3}\right\rceil,\ldots,2n\big\}} P^o_{\omega} (X_{\ell} = b_+) \label{I-eq2}
\end{align}
where we used \eqref{I-eq1} in the third step. Here, for $x,y \in \Z$,
\[
\widehat{\inf}_{\ell \in \big\{\left\lceil \tfrac{4n}{3} \right\rceil,\ldots,2n\big\}} P^x_{\omega}(X_{\ell}=y)
\] 
is the short notation for 
\[
\inf_{\ell \in \big\{\left\lceil \tfrac{4n}{3}\right\rceil,\ldots,2n\big\}\cap \big(2\Z+ (x+y)\big)} P^x_{\omega}(X_{\ell}=y)
\]
since we have to take care of the periodicity of the random walk.
\\
Let us now have a closer look at the factors in the lower bound in \eqref{I-eq2} separately:\\[10pt]
\underline{First factor in \eqref{I-eq2}:}
\\We can bound the first factor from below by
\begin{align*}
& P^o_{\omega} \left(\tau(b_+) \le \tfrac{2n}{3},\ \tau(b_+) < \tau(b_-) \right) \vphantom{\exp \left( \max_{b_- \le i \le j \le b_+} \big(V(j) - V(i)\big)  \right)}\\
\ge\ & 1 - P^o_{\omega} \left(\tau(b_+) > \tfrac{2n}{3},\ \tau(b_+) < \tau(b_-) \right) - P^o_{\omega} \big( \tau(b_+) \ge \tau(b_-) \big)\vphantom{\exp \left( \max_{b_- \le i \le j \le b_+} \big(V(j) - V(i)\big)  \right)}\\
\ge\ & 1 - \tfrac{3}{2n}\cdot E^o_{\omega} \left[\tau(b_+) \cdot \mathds{1}_{\{\tau(b_+) < \tau(b_-)\}} \right] - P^o_{\omega} \big( \tau(b_+) \ge \tau(b_-) \big)\vphantom{\exp \left( \max_{b_- \le i \le j \le b_+} \big(V(j) - V(i)\big)  \right)}\\
\ge\ & 1 - \tfrac{3}{2n}\cdot (b_+ - b_-)^2\cdot \exp \left( \max_{b_- \le i \le j \le b_+} \big(V(j) - V(i)\big)  \right) - \frac12 \ ,
\end{align*}
where we used \eqref{I-prel4} and assumption \eqref{I-ass4} for the last step.
Therefore, we get for \mbox{$\omega \in \Gamma(L,\delta)$} and $\exp\left(3 \delta L \right) \le n$ that
\begin{align}
&  P^o_{\omega} \left(\tau(b_+) \le \tfrac{2n}{3},\ \tau(b_+) < \tau(b_-) \right)  \ge \frac12 - \frac{3 \cdot 4 \cdot L^4}{2\cdot \exp(3\delta L)} \cdot \exp(2\delta L)= \frac12 - 6 \cdot L^4 \cdot \exp(-\delta L). \label{I-eq5}
\end{align}
\underline{Second factor in \eqref{I-eq2}:}
\\By using assumption \eqref{I-ass2} and the relation in \eqref{I-eq1.1}, we get for $\omega \in \Gamma(L,\delta)$:
\begin{align}
&\frac{\mu_{\omega}(0)}{\mu_{\omega}(b_+)} =  \frac{\tfrac{1}{\omega_0}}{e^{-V(b_+)} + e^{-V(b_+-1)}} = \frac{\tfrac{1}{\omega_0}}{e^{-V(b_+)}\cdot (1+\rho_{b_+})} \nonumber \\
\ge\ & \frac{\tfrac{1}{1-\varepsilon}}{1+\tfrac{1-\varepsilon}{\varepsilon}} \cdot e^{V(b_+)}= \frac{\varepsilon}{1 - \varepsilon} \cdot e^{V(b_+)} \ge \frac{\varepsilon}{1 - \varepsilon} \cdot \exp(-\delta L). \label{I-eq3}
\end{align}
Here we used that $V(b_+) \ge -\delta L$ holds for $\omega \in \Gamma(L,\delta)$.\\[10pt]
\underline{Third factor in \eqref{I-eq2}:}\\
For the last factor in \eqref{I-eq2}, we can compare the RWRE with the process $(\widetilde{X}_n)_{n \in \N_0}$ which behaves as the original RWRE but is reflected at the positions $T^-:=T^-(L)$ and $T^+:=~T^+(L)$, i.e.\ we have for $x \in \{T^-,T^-+1,\ldots,T^+\}$
\begin{align*}
& P_{\omega}^{x}(\widetilde{X}_0=x) = 1 \vphantom{\widehat{\inf}_{k \in \big\{\lceil \tfrac{\ell}{2}\rceil,\ldots,\ell\big\}}},\\
& P_{\omega}^{x}(\widetilde{X}_{n+1} = y\pm1|\widetilde{X}_n=y) = P_{\omega}^{x}(X_{n+1} = y\pm1|X_n=y)  \vphantom{\widehat{\inf}_{k \in \big\{\lceil \tfrac{\ell}{2}\rceil,\ldots,\ell\big\}}},\\
& & \hspace{-13.98193pt}  \forall y \in \{T^-+1,T^-+2,\ldots,T^+-1\} \vphantom{\widehat{\inf}_{k \in \big\{\lceil \tfrac{\ell}{2}\rceil,\ldots,\ell\big\}}},\\
& P_{\omega}^{x}(\widetilde{X}_{n+1} = y+1|\widetilde{X}_n=y) = 1 \quad & \text{for } y= T^- \vphantom{\widehat{\inf}_{k \in \big\{\lceil \tfrac{\ell}{2}\rceil,\ldots,\ell\big\}}},\\
& P_{\omega}^{x}(\widetilde{X}_{n+1} = y-1|\widetilde{X}_n=y) = 1 \quad  & \text{for } y= T^+. \vphantom{\widehat{\inf}_{k \in \big\{\lceil \tfrac{\ell}{2}\rceil,\ldots,\ell\big\}}}
\end{align*}
Therefore, we have for $\ell \in \big\{\left\lceil \tfrac{4n}{3}\right\rceil,\ldots,2n\big\}\cap \big(2\Z+ b_+\big)$
\begin{align}
&P^{o}_{\omega} (X_{\ell} = b_+) \nonumber \vphantom{\widehat{\inf}_{k \in \big\{\lceil \tfrac{\ell}{2}\rceil,\ldots,\ell\big\}}}\\
\ge\ &  P^{o}_{\omega} (X_{\ell} = b_+,\ \min\{\tau(T^-),\tau(T^+)\} > 2n) \nonumber \vphantom{\widehat{\inf}_{k \in \big\{\lceil \tfrac{\ell}{2}\rceil,\ldots,\ell\big\}}} \\
=\ & P^{o}_{\omega} (\widetilde{X}_{\ell} = b_+) - P^{o}_{\omega} (\widetilde{X}_{\ell} = b_+,\ \min\{\tau(T^-),\tau(T^+)\} \le 2n) \nonumber \vphantom{\widehat{\inf}_{k \in \big\{\lceil \tfrac{\ell}{2}\rceil,\ldots,\ell\big\}}}\\
\ge\ & P^{o}_{\omega} (\widetilde{X}_{\ell} = b_+) - P^{o}_{\omega} ( \min\{\tau(T^-),\tau(T^+)\} \le 2n) \nonumber \vphantom{\widehat{\inf}_{k \in \big\{\lceil \tfrac{\ell}{2}\rceil,\ldots,\ell\big\}}}\\
\ge\ & P^{o}_{\omega} \big(\widetilde{X}_{\ell} = b_+, \ \tau(b_+) \le \tfrac{\ell}{2},\ \tau(b_+) < \tau(b_-)\big) - P^{o}_{\omega} ( \min\{\tau(T^-),\tau(T^+)\} \le 2n) \nonumber \vphantom{\widehat{\inf}_{k \in \big\{\lceil \tfrac{\ell}{2}\rceil,\ldots,\ell\big\}}}\\
\ge\ & P^{o}_{\omega} \big(\tau(b_+) \le \tfrac{\ell}{2},\ \tau(b_+) < \tau(b_-)\big) \cdot \widehat{\inf\limits}_{k \in \big\{\left\lceil \tfrac{\ell}{2}\right\rceil,\ldots,\ell\big\}} P^{b_+}_{\omega} (X_{k} = b_+) \nonumber \vphantom{\widehat{\inf}_{k \in \big\{\lceil \tfrac{\ell}{2}\rceil,\ldots,\ell\big\}}}\\
& - P^{o}_{\omega} ( \min\{\tau(T^-),\tau(T^+)\} \le 2n).  \label{I-eq4} \vphantom{\widehat{\inf}_{k \in \big\{\lceil \tfrac{\ell}{2}\rceil,\ldots,\ell\big\}}}
\end{align}
Using \eqref{I-prel2} and \eqref{I-prel3}, we see that the last term in \eqref{I-eq4} with the negative sign decreases exponentially for $n \le \exp\big((1-2\delta)L\big)$, i.e.
\begin{align}
& P^{o}_{\omega} ( \min\{\tau(T^-),\tau(T^+)\} \le 2n) \le P^{o}_{\omega} \left( \min\{\tau(T^-),\tau(T^+)\} \le 2\cdot\exp\big((1-2\delta)L\big) \right) \nonumber \vphantom{\widehat{\inf}_{k \in \big\{\lceil \tfrac{\ell}{2}\rceil,\ldots,\ell\big\}}}\\
\le\ & P^{o}_{\omega} \left( \tau(T^-) \le 2\cdot\exp\big((1-2\delta)L\big) \right) + P^{o}_{\omega} \left( \tau(T^+) \le 2\cdot\exp\big((1-2\delta)L\big)\right) \nonumber \vphantom{\widehat{\inf}_{k \in \big\{\lceil \tfrac{\ell}{2}\rceil,\ldots,\ell\big\}}}\\
\le\ & 4 \cdot \exp\big((1-2\delta)L\big) \cdot \exp\big(-L\big) = 4 \cdot \exp\big(-2\delta L\big) \label{I-eq4a}. \vphantom{\widehat{\inf}_{k \in \big\{\lceil \tfrac{\ell}{2}\rceil,\ldots,\ell\big\}}}
\end{align}
In order to derive a lower bound for the first term in \eqref{I-eq4}, we first notice that the analogous calculation as in \eqref{I-eq5} shows for $\omega \in \Gamma(L,\delta)$ that
\begin{align}
&  P^o_{\omega} \left(\tau(b_+) \le \tfrac{\ell}{2},\ \tau(b_+) < \tau(b_-) \right) \ge 1 - \frac{2}{\ell} \cdot 4 \cdot L^4 \cdot \exp(2\delta L) - \frac12 \nonumber\\
\ge\ &  \frac12 - 6 \cdot L^4 \cdot \exp(-\delta L) \label{I-eq5.1}
\end{align}
since $\ell \ge \left\lceil \tfrac{4n}{3}\right\rceil \ge \tfrac43 \cdot \exp(3 \delta L)$ for $n \ge \exp(3 \delta L)$. For the second factor, we show the following\\

\begin{lem} \label{I-lem3}
For $\omega \in \Gamma(L,\delta)$ and for all $\ell \in 2 \N$, we have
\[
P^{b_+}_{\omega}(\widetilde{X}_{\ell} = b_+) \ge \frac12 \cdot \frac{1}{|T^-|+T^+ +1} \cdot \exp\big(-\delta L\big).
\]
\end{lem}
\begin{pfof}{Lemma \ref{I-lem3}}
Using the reversibility (cf.\ \eqref{I-eq1}) of $(\widetilde{X}_{\ell})_{\ell \in \N_0}$, we get
\begin{align}
& P^{b_+}_{\omega}(\widetilde{X}_{\ell} = b_+) \nonumber \vphantom{\sum_{x=T^-}^{T^+}} \allowdisplaybreaks[0]\\
=\ & \sum_{x=T^-}^{T^+} P^{b_+}_{\omega}(\widetilde{X}_{\ell/2} = x) \cdot P^x_{\omega}(\widetilde{X}_{\ell/2} = b_+) \nonumber\\
=\ & \sum_{x=T^-}^{T^+} P^{b_+}_{\omega}(\widetilde{X}_{\ell/2} = x) \cdot \frac{\widetilde{\mu}_{\omega}(b_+)}{\widetilde{\mu}_{\omega}(x)} \cdot P^{b_+}_{\omega}(\widetilde{X}_{\ell/2} = x), \label{I-eq7}
\end{align}
where $\widetilde{\mu}_{\omega}(\cdot)$ denotes a reversible measure of the reflected random walk $(\widetilde{X}_n)_{n \in \N_0}$ which is unique up to multiplication by a constant. To see that $(\widetilde{X}_{\ell})_{\ell \in \N_0}$ is also reversible, it is enough to note that $(\widetilde{X}_{\ell})_{\ell \in \N_0}$ can again be described as an electrical network with the following conductances:
\begin{align*}
\widetilde{C}_{(x,x+1)}(\omega)= \begin{cases}
C_{(x,x+1)}(\omega) = e^{-V(x)} & \text{for } x=T^-, T^-+1,\ldots, T^+-1 \\
0 & \text{for } x = T^- - 1, T^+
\end{cases}
\end{align*}
Therefore, a reversible measure for the reflected random walk is given by (cf.\ \eqref{I-eq1.1})
\[
\widetilde{\mu}_{\omega}(x) = \begin{cases}\mu_{\omega}(x) = e^{-V(x)} + e^{-V(x-1)} & \text{for } x=T^-+1, T^-+2,\ldots, T^+-1, \\
e^{-V(T^-)} & \text{for } x=T^-, \\
e^{-V(T^+-1)} & \text{for } x=T^+.
\end{cases}
\]
This implies, since $0 \le b_+ < T^+$,
\begin{align}
& \frac{\widetilde{\mu}_{\omega}(b_+)}{\widetilde{\mu}_{\omega}(x)} \ge \frac{e^{-V(b_+)}+ e^{-V(b_+-1)}}{e^{-V(x)}+ e^{-V(x-1)}} \nonumber \\
\ge\ & \frac{e^{-V(b_+)}}{2\cdot e^{\left( - \min \{V(b_+),V(b_-) \} \right)}} \ge \frac{1}{2} \cdot \exp(-\delta L) \label{I-eq8}
\end{align}
for $T^- \le x \le T^+$ and for $\omega \in \Gamma(L,\delta)$. By applying \eqref{I-eq8} to \eqref{I-eq7}, we get
\begin{align}
& P^{b_+}_{\omega}(\widetilde{X}_{\ell} = b_+) \nonumber \\
\ge\ &\frac12 \cdot \sum_{x=T^-}^{T^+} \left(P^{b_+}_{\omega}(\widetilde{X}_{\ell/2} = x)\right)^2 \cdot \exp(-\delta L)\nonumber\\
\ge\ &\frac12 \cdot \sum_{x=T^-}^{T^+} \left(\frac{1}{|T^-|+T^+ +1}\right)^2 \cdot \exp(-\delta L) \nonumber\\
=\ & \frac12 \cdot \frac{1}{|T^-|+T^+ +1} \cdot \exp(-\delta L). \label{I-eq9}
\end{align}
Here, we used that we have
\begin{align*}
\sum_{x=T^-}^{T^+} (a_x)^2 \ge \sum_{x=T^-}^{T^+} \left(\frac{1}{|T^-|+T^+ +1}\right)^2
\end{align*}
for every sequence $(a_x)_{x}$ with $\sum\limits_{x=T^-}^{T^+} a_x =1$.
\renewcommand{\qedsymbol}{\hfill$\square$\vspace{1ex}}
\end{pfof}

We can now return to the proof of Proposition \ref{I-prop1} and finish our lower bound for the third factor in \eqref{I-eq2}. By applying \eqref{I-eq4a}, \eqref{I-eq5.1} and Lemma \ref{I-lem3} to \eqref{I-eq4}, we get for $\exp(3\delta L) \le n \le \exp((1-2\delta) L)$ and $\omega \in \Gamma(L,\delta)$, i.e.\ $|T^-|, T^+ \le L^2$,
\begin{align}
& \widehat{\inf}_{ \ell \in \big\{\left\lceil \tfrac{4n}{3}\right\rceil,\ldots,2n\big\}}P^{o}_{\omega} (X_{\ell} = b_+) \nonumber \vphantom{\frac12}\\
\ge\ &  \left(\frac12 - 6 \cdot L^4 \cdot \exp(-\delta L)\right) \cdot \frac12 \cdot \frac{1}{2L^2 +1} \cdot \exp(-\delta L) - 4 \cdot \exp\big(-2\delta L\big) \nonumber \\
\ge\ &  \exp\left(-\tfrac32\delta L\right)\label{I-eq12} \vphantom{\frac12}
\end{align}
for all $L=L(\delta)$ large enough.\\[10pt]
To finish the proof of Proposition \ref{I-prop1}, we can collect our lower bounds in \eqref{I-eq5}, \eqref{I-eq3}, and \eqref{I-eq12} and conclude with \eqref{I-eq2} that for $\exp\left(3 \delta L \right) \le n \le \exp\big((1-2\delta)L\big)$ and for $\omega \in \Gamma(L,\delta)$ we have
\begin{align*}
& P_{\omega} (X_{2n}=0) \vphantom{\frac12}\\
\ge\ & \left(\frac12 - 6 \cdot L^4 \exp(-\delta L)\right) \cdot \frac{\varepsilon}{1-\varepsilon} \exp(-\delta L) \cdot \exp\left(-\tfrac32\delta L\right)\\
\ge\ &  \exp(-3\delta L)  \vphantom{\frac12}
\end{align*}
for all $L=L(\delta)$ large enough. This shows \eqref{I-lem1} since we have $P_{\omega}(X_{2n}=0) \ge \varepsilon^{2n} > 0$ for all $n \in \N$ due to assumption \eqref{I-ass2}.
\end{pfof}
\begin{prop} \label{I-lemma2}
For $0 < \delta < 1$, we have
\begin{equation}
\label{I-lem2}
\p (\omega:\ \omega \in \Gamma(L,\delta) \text{ for infinitely many }L )=1.
\end{equation}
\end{prop}
\begin{pfof}{Proposition \ref{I-lemma2}}
Let $(B(t))_{t \in \R}$ be the two-sided Brownian motion from Theorem \ref{I-Komlos} and let us choose some $0 < \delta < \tfrac12$. For $y \in \R$ we define
\begin{align*}
&\widehat{T}^+(y):= \inf\{t \ge 0:\ B(t)= y\},\\
&\widehat{T}^-(y):= \sup\{t \le 0:\ B(t)= y\}
\end{align*}
as the first hitting times of $y$ on the positive and negative side of the origin, respectively. Additionally, for $L \in \N$, $i \in \N$, $y \in \R$, we can introduce the following sets
\begin{align*}
F_L^+(y):=\ & \{\widehat{T}^+ \left(y \cdot L \right) < \widehat{T}^+ \left(- y \cdot L \right)\}, \\
F_L^-(y):=\ & \{\widehat{T}^-\left(y \cdot L \right) < \widehat{T}^-\left(-y \cdot L \right)\}
\end{align*}
on which the Brownian motion reaches the value $y \cdot L$ before $-y \cdot L$. Further we define
\begin{align*}
G_L^+(i):=\ & \left\{B(t) \ge (2i-1) \cdot \tfrac{\delta}{4} \cdot  L  \quad \text{for} \quad \widehat{T}^+\left(2i \cdot \tfrac{\delta}{4} \cdot L \right) \le t \le \widehat{T}^+\left((2i+2) \cdot \tfrac{\delta}{4} \cdot L \right) \right\}, \\
G_L^-(i):=\ & \left\{B(t) \ge (2i-1) \cdot \tfrac{\delta}{4} \cdot  L  \quad \text{for} \quad \widehat{T}^-\left((2i+2) \cdot \tfrac{\delta}{4} \cdot L \right) \le t \le \widehat{T}^-\left(2i \cdot \tfrac{\delta}{4} \cdot L \right) \right\}
\intertext{on which the Brownian motion does not decrease much between the first hitting time of the two levels of interest. Using these sets, we can define the sets}
A^+(L,\delta):=\ & F_L^+(\delta) \cap \left\{ \widehat{T}^+(1.1 \cdot L) \le L^2,\  \min_{\widehat{T}^+(\delta \cdot L) \le t \le \widehat{T}^+(1.1 \cdot L)} B(t) \ge \frac{\delta}{4} \cdot  L \right\}, \\
A^-(L,\delta):=\ & F_L^-(\delta) \cap \left\{ - \widehat{T}^-(1.1 \cdot L) \le L^2,\  \min_{\widehat{T}^-(1.1 \cdot L) \le t \le \widehat{T}^-(\delta \cdot L)} B(t) \ge \frac{\delta}{4} \cdot  L \right\}, \\
D^+(L,\delta):=\ & G_L^+(0) \cap G_L^+(1) \cap G_L^+(2) \allowdisplaybreaks[0]\\
& \cap \left\{ \widehat{T}^+(1.2 \cdot L) \le 0.9 \cdot L^2,\  \min_{\widehat{T}^+ \big(\tfrac{3 \cdot \delta}{2} \cdot L \big) \le t \le \widehat{T}^+(1.2 \cdot L)} B(t) \ge \frac{3\delta}{4} \cdot  L \right\}, \\
D^-(L,\delta):=\ & G_L^-(0) \cap G_L^-(1) \cap G_L^-(2) \\
& \cap \left\{ - \widehat{T}^-(1.2 \cdot L) \le 0.9 \cdot L^2,\  \min_{\widehat{T}^-(1.2 \cdot L) \le t \le \widehat{T}^-\big(\tfrac{3\delta}{2} \cdot L \big)} B(t) \ge \frac{3\delta}{4} \cdot  L \right\}
\end{align*}
which which will be used for an approximation of our previously constructed valleys $\omega$ belonging to $\Gamma(L,\delta)$ which we illustrated in Figure \ref{I-figure1} on page \pageref{I-figure1}. Here, we added the factors $1.1$, $1.2$ and $0.9$ in contrast to the construction before in order to have some space for the approximation. For the Brownian motion, we can directly compute that we have
\begin{equation}
\label{I-eq15}
\p\big(D^+(1,\delta) \cap D^-(1,\delta)\big) > 0.
\end{equation}
Thereby, the scaling property of the Brownian motion, i.e. the property that for $L \in \N$
\begin{equation} \label{I-eq16a}
\left( \frac{1}{L} B(L^2 \cdot t) \right)_{ t \in \R}
\end{equation}
is again a two-sided Brownian motion with diffusion constant $\sigma$, implies 
\begin{equation} \label{I-eq16}
\p\big(D^+(L,\delta) \cap D^-(L,\delta)\big) = \p\big(D^+(1,\delta) \cap D^-(1,\delta)\big) > 0
\end{equation}
for all $L \in \N$. 

At first, we notice that for $L_0 \in \N$ we have
\begin{align} \label{I-eq25}
& \p \left(\, \bigcap_{L=L_0}^{\infty} \Big( A^+(L,\delta) \cap A^-(L,\delta) \Big)^c \right)
\le \p \left(\ \bigcap_{k=\ell+1}^{\infty} \Big( A^+(L_k,\delta) \cap A^-(L_k,\delta) \Big)^c \right)
\end{align} 
for arbitrary $\ell \in \N_0$, where we define 
\[
L_k := \max \left\{10, \left\lceil \tfrac{2}{\delta} \right\rceil \right\} \cdot (L_{k-1})^2
\]
for $k \in \N$ inductively. Note that for $n > \ell +1$ 
with
\[
\mathcal{F}_n:=\sigma \left( \big(B(t)\big)_{-(L_{n-1})^2 \le t \le (L_{n-1})^2} \right),
\]
the following holds:
\begin{align}
& \p \left(\ \bigcap_{k=\ell+1}^{n} \Big( A^+(L_k,\delta) \cap A^-(L_k,\delta) \Big)^c \right) \nonumber \vphantom{\E \left[ \prod_{k=\ell+1}^{n-1} \mathds{1}_{\big( A^+(L_k,\delta) \cap A^-(L_k,\delta)  \big)^c} \cdot \mathds{1}_{\left\{ \max\limits_{-(L_{n-1})^2 \le t \le (L_{n-1})^2} |B(t)| \le (L_{n-1})^2 \right\}} \right.}\\
\le\ &  \E \left[ \prod_{k=\ell+1}^{n-1} \mathds{1}_{\big( A^+(L_k,\delta) \cap A^-(L_k,\delta)  \big)^c} \cdot \mathds{1}_{\left\{ \max\limits_{-(L_{n-1})^2 \le t \le (L_{n-1})^2} |B(t)| < (L_{n-1})^2 \right\}} \right. \nonumber\allowdisplaybreaks[0]\\
& \hspace{-5.5pt} \cdot \left. \left.  \E\left[ \vphantom{\prod_{k=\ell+1}^{n-1}} \mathds{1}_{\left\{ \big(B(t + (L_{n-1})^2) - B((L_{n-1})^2) \big)_{t \in \R} \notin D^+ \left(L_n, \delta \right) \right\} \cup \left\{ \big(B(t - (L_{n-1})^2) - B(-(L_{n-1})^2) \big)_{t \in \R} \notin D^- \left(L_n, \delta \right) \right\}}  \right| \hspace{-1pt} \mathcal{F}_n \hspace{-1pt} \right] \hspace{-2pt} \right] \nonumber\allowdisplaybreaks[0]\\
& +  \p \left( \max\limits_{-(L_{n-1})^2 \le t \le (L_{n-1})^2} |B(t)| \ge (L_{n-1})^2 \right) \nonumber \vphantom{\E \left[ \prod_{k=\ell+1}^{n-1} \mathds{1}_{\big( A^+(L_k,\delta) \cap A^-(L_k,\delta)  \big)^c} \cdot \mathds{1}_{\left\{ \max\limits_{-(L_{n-1})^2 \le t \le (L_{n-1})^2} |B(t)| \le (L_{n-1})^2 \right\}} \right.}\\
\le\ &  \Big(1 - \p \Big(D^+ \left(L_n , \delta \right) \cap D^- \left(L_n, \delta \right) \Big) \Big) \cdot \p \left(\ \bigcap_{k=\ell+1}^{n-1} \Big( A^+(L_k,\delta) \cap A^-(L_k,\delta) \Big)^c \right) \nonumber \vphantom{\E \left[ \prod_{k=\ell+1}^{n-1} \mathds{1}_{\big( A^+(L_k,\delta) \cap A^-(L_k,\delta)  \big)^c} \cdot \mathds{1}_{\left\{ \max\limits_{-(L_{n-1})^2 \le t \le (L_{n-1})^2} |B(t)| \le (L_{n-1})^2 \right\}} \right.}\\
&  + \p \left( \max\limits_{-(L_{n-1})^2 \le t \le (L_{n-1})^2} |B(t)| \ge (L_{n-1})^2 \right) \nonumber \vphantom{\E \left[ \prod_{k=\ell+1}^{n-1} \mathds{1}_{\big( A^+(L_k,\delta) \cap A^-(L_k,\delta)  \big)^c} \cdot \mathds{1}_{\left\{ \max\limits_{-(L_{n-1})^2 \le t \le (L_{n-1})^2} |B(t)| \le (L_{n-1})^2 \right\}} \right.}\\
\le\ &  \Big(1 - \p \Big(D^+ \left(1 , \delta \right) \cap D^- \left(1, \delta \right) \Big) \Big)^{n-\ell}   + \sum_{k=\ell+1}^{n} \p \left( \max\limits_{-(L_{k-1})^2 \le t \le (L_{k-1})^2} |B(t)| \ge (L_{k-1})^2 \right) \label{I-eq26} \vphantom{\E \left[ \prod_{k=\ell+1}^{n-1} \mathds{1}_{\big( A^+(L_k,\delta) \cap A^-(L_k,\delta)  \big)^c} \cdot \mathds{1}_{\left\{ \max\limits_{-(L_{n-1})^2 \le t \le (L_{n-1})^2} |B(t)| \le (L_{n-1})^2 \right\}} \right.}.
\end{align}
To see that the first step holds, note that for 
\begin{align}
\omega \in &
\left\{\max\limits_{-(L_{n-1})^2 \le t \le (L_{n-1})^2)} |B(t)| < (L_{n-1})^2 \right\} \nonumber\\
& \cap \left\{\big(B(t + (L_{n-1})^2 ) - B((L_{n-1})^2) \big)_{t \in \R} \in D^+ \left(L_n, \delta \right) \right\} \label{I-eq26.x}
\end{align}
we have
\begin{align*}
 \min_{0 \le t \le (L_n)^2} B(t) \ge\ & \min_{0 \le t \le (L_{n-1})^2} B(t) + \min_{(L_{n-1})^2 \le t \le (L_n)^2} B(t+(L_{n-1})^2) - B((L_{n-1})^2) \\
>\ &  - (L_{n-1})^2 - \frac{\delta}{4} \cdot L_n > - \delta \cdot L_n, 
\intertext{and}
\max_{0 \le t \le (L_n)^2} B(t) \ge\ & B\big((L_{n-1})^2\big) + \max_{(L_{n-1})^2 \le t \le (L_n)^2 - (L_{n-1})^2} B(t+(L_{n-1})^2) - B((L_{n-1})^2)\\
\ge\ &  - (L_{n-1})^2 + 1.2 \cdot L_n \ge 1.1 \cdot L_n.
\end{align*}
In particular, we have $\widehat{T}^+(\delta \cdot L) < \widehat{T}^+(- \delta \cdot L)$ and $\widehat{T}^+(1.1 \cdot L) \le L^2$ on the considered set. Similarly, again on the set in \eqref{I-eq26.x}, we see that we have
\begin{align*}
&  \widehat{T}^+(\delta \cdot L) > \inf\{t \ge (L_{n-1})^2:\ (B(t + (L_{n-1})^2 ) - B((L_{n-1})^2) \ge \tfrac{\delta}{2} \cdot L_n\}, \\
&  \widehat{T}^+(\delta \cdot L) < \inf\{t \ge (L_{n-1})^2:\ (B(t + (L_{n-1})^2 ) - B((L_{n-1})^2) \ge \tfrac{3\cdot \delta}{2} \cdot L_n\}, 
\intertext{which implies}
&  \min_{\widehat{T}^+(\delta \cdot L) \le t \le \widehat{T}^+(1.1 \cdot L)} B(t) \ge  \frac{\delta}{4} \cdot L_n
\end{align*}
by construction of $D^+(L_n,\delta)$.
Altogether, we can conclude that $\omega \in A^+(L_n,\delta)$ holds for our choice of $\omega$ in \eqref{I-eq26.x}. The argument for the negative part runs completely analogously. Further in \eqref{I-eq26}, we used the Markov property of the Brownian motion in the second step. Additionally, we iterated the first two steps $n-\ell-1$ times and used \eqref{I-eq16} for the last step. To control the last sum in \eqref{I-eq26}, let us recall the standard upper bound 
\[
\p\left( Z \ge x \right) \le \frac{1}{x} \cdot \frac{1}{\sqrt{2\pi}} \cdot \exp \left(- \frac{x^2}{2} \right) \quad \text{for } x > 0
\]
for a random variable $Z \sim \mathcal{N}(0,1)$, which can be found for example in Lemma 12.9 in Appendix B of \cite{mörters}. By using this upper bound, we can conclude that
\begin{align}
& \sum_{k=\ell+1}^{n} \p \left( \max\limits_{-(L_{k-1})^2 \le t \le (L_{k-1})^2} |B(t)| \ge (L_{k-1})^2 \right) \le 4 \cdot \sum_{k=\ell+1}^{n} \p \left( \max\limits_{0\le t \le (L_{k-1})^2} \frac{B(t)}{\sigma \cdot L_{k-1}}  \ge \frac{L_{k-1}}{\sigma} \right) \nonumber \\
\le\ & 4 \cdot \sum_{k=\ell+1}^{\infty} \frac{\sigma}{L_{k-1}} \cdot \frac{1}{\sqrt{2 \pi}} \cdot \exp \left(- \frac{(L_{k-1})^2}{2 \sigma^2}  \right) \xrightarrow{\ell \to \infty} 0. \label{I-eq27}
\end{align}
Here, we used that 
\[
\max\limits_{0\le t \le (L_{k-1})^2} \frac{B(t)}{\sigma \cdot L_{k-1}} \sim |Z|
\]
for all $k \in \N$, where $Z \sim \mathcal{N}(0,1)$. By combining the upper bounds in \eqref{I-eq25}, \eqref{I-eq26}, and \eqref{I-eq27}, we get for all $\ell \in \N_0$
\begin{align*}
& \p \left( \omega \notin  \big( A^+(L,\delta) \cap A^-(L,\delta) \big) \text{ for all } L \ge L_0\right) \vphantom{\sum_{k=\ell+1}^{\infty}}\\
\le\ & \lim_{n \to \infty}  \Big(1 - \p \Big(D^+ \left(1 , \tfrac{\delta}{2}\right) \cap D^- \left(1,\tfrac {\delta}{2} \right) \Big) \Big)^{n-\ell} \vphantom{\sum_{k=\ell+1}^{\infty}} \\
&  + \sum_{k=\ell+1}^{\infty} \p \left( \max\limits_{-(L_{k-1})^2 \le t \le (L_{k-1})^2} |B(t)| \ge (L_{k-1})^2 \right) \xrightarrow{\ell \to \infty} 0.
\end{align*}
Since $L_0 \in \N$ was chosen arbitrarily, we can conclude that for $0 < \delta < \tfrac12$ we have
\[
\p \left(\omega:\  \omega \in  \big(A^+(L, \delta) \cap A^-(L, \delta)\big) \text{ for infinitely many } L \right) = 1.
\]
Using the Koml{\'o}s-Major-Tusn{\'a}dy strong approximation Theorem (cf.\ Theorem \ref{I-Komlos}), we see that for $0 < \delta < \tfrac12$ we have
\begin{align*}
& \left\{\omega:\  \omega \in  \big(A^+(L, \delta) \cap A^-(L, \delta)\big) \text{ for infinitely many } L \right\} \allowdisplaybreaks[0]\\
\subseteq\ &  \left\{\omega:\ \omega \in \Gamma(L,2\delta) \text{ for infinitely many }L  \right\}, 
\end{align*}
which is enough to conlude that \eqref{I-lem2} holds for all $0 < \delta < 1$.
\end{pfof}

With the help of Proposition \ref{I-prop1} and Proposition \ref{I-lemma2}, we can now turn to the proofs of our Theorems \ref{I-Rthm1} --  \ref{I-Rthm3}:

\begin{pfof}{Theorem \ref{I-Rthm1}} For a fixed $0 \le \alpha < 1$, we choose $0 < \delta < \tfrac16$ such that
\begin{equation}
\label{I-prthm1}
\alpha < \frac{1-5\delta}{1-2\delta}\ . 
\end{equation}
For $\omega \in \Gamma(L, \delta)$, the inequality in \eqref{I-lem1} implies that
\begin{align}
& \hphantom{\ge}\ \sum_{n \in \N} P_{\omega}(X_{2n}=0) \cdot n^{-\alpha} \ge \sum_{\lceil\exp(3\delta L)\rceil \le n \le \lfloor\exp((1-2\delta)L)\rfloor} P_{\omega}(X_{2n}=0) \cdot n^{-\alpha} \nonumber\\
&\ge \Big(\exp\big((1-2\delta)L\big) - \exp(3\delta L) - 1 \Big) \cdot C \cdot \exp(-3\delta L) \cdot \left(\exp\big((1-2\delta)L\big)\right)^{-\alpha} \nonumber \vphantom{\sum_{\lceil\exp(3\delta L)\rceil \le n \le \lfloor\exp((1-2\delta)L)\rfloor}}\\
&= C \cdot \exp(-3\delta L) \cdot \exp(3\delta L) \cdot \Big(\exp\big((1-5\delta)L\big) - 1 - \exp(-3\delta L)\Big) \cdot \exp\big(-\alpha(1-2\delta)L\big) \nonumber \vphantom{\sum_{\lceil\exp(3\delta L)\rceil \le n \le \lfloor\exp((1-2\delta)L)\rfloor}}\allowdisplaybreaks[0]\\
& \xrightarrow{L\to\infty} \infty \nonumber \vphantom{\sum_{\lceil\exp(3\delta L)\rceil \le n \le \lfloor\exp((1-2\delta)L)\rfloor}}.
\end{align}
Since Proposition \ref{I-lemma2} shows that for $\p$-a.e.\ environment $\omega$ we find $L$ arbitrarily large such that $\omega \in \Gamma(L,\delta)$, we can conclude that \eqref{I-thm1} holds for $\p$-a.e.\ environment $\omega$.
\end{pfof}

\begin{pfof}{Theorem \ref{I-Rthm2}}
For fixed $\alpha > 0$, we choose $\delta$ such that
\begin{align*}
0 < \delta < \min \left\{ \frac{1}{2+3\alpha}, \frac{1}{5} \right\}, \intertext{which yields} 1-2\delta -3\alpha \delta > 0 \quad  \text{and} \quad  1- 2 \delta > 3 \delta.
\end{align*}
For $\omega \in \Gamma(L,\delta)$, the inequality in \eqref{I-lem1} implies
\begin{align*}
& \hphantom{\ge}\ \sum_{n \in \N} \Big( P_{\omega}(X_{2n}=0)\Big)^{\alpha} \ge \sum_{\lceil\exp(3\delta L)\rceil \le n \le \lfloor\exp((1-2\delta)L)\rfloor} \Big(P_{\omega}(X_{2n}=0)\Big)^{\alpha} \\
&\ge \Big(\exp\big((1-2\delta)L\big) - \exp(3\delta L) - 1 \Big) \cdot \big(C \cdot \exp(-3\delta L)\big)^{\alpha} \vphantom{\sum_{\lceil\exp(3\delta L)\rceil \le n \le \lfloor\exp((1-2\delta)L)\rfloor}}\\
&=  C^{\alpha} \cdot  \exp(-3 \alpha \delta L) \cdot\exp(3 \alpha \delta L)   \vphantom{\sum_{\lceil\exp(3\delta L)\rceil \le n \le \lfloor\exp((1-2\delta)L)\rfloor}}\\
& \qquad \cdot \Big(\exp\big((1-2\delta -3\alpha \delta)L\big) - \exp\big((3\delta - 3 \alpha \delta) L\big)  -  \exp( -3\alpha \delta L) \Big) \vphantom{\sum_{\lceil\exp(3\delta L)\rceil \le n \le \lfloor\exp((1-2\delta)L)\rfloor}}\\
&\xrightarrow{L\to\infty} \infty. \vphantom{\sum_{\lceil\exp(3\delta L)\rceil \le n \le \lfloor\exp((1-2\delta)L)\rfloor}}
\end{align*}
Again since Proposition \ref{I-lemma2} shows that for $\p$-a.e.\ environment $\omega$ we find $L$ arbitrarily large such that $\omega \in \Gamma(L,\delta)$, we can conclude that \eqref{I-thm2} holds for $\p$-a.e.\ environment $\omega$.
\end{pfof}

\begin{pfof}{Theorem \ref{I-Rthm3}}
Due to the independence of the environments $\omega^{(1)}, \omega^{(2)}, \ldots, \omega^{(d)}$, we can extend the proof of Proposition \ref{I-lemma2} to get
\begin{equation}
\label{I-eq16.2}
\p^{\otimes d} \left(\text{For infinitely many } L \in \N \text{ we have}:\ \omega^{(i)} \in \Gamma(L,\delta)\ \text{for } i=1,2,\ldots d \right)=1
\end{equation}
for all $0 < \delta < 1$.\\[10pt]
Thereby due to Proposition \ref{I-prop1}, we have for $(\omega^{(1)}, \omega^{(2)}, \ldots, \omega^{(d)})$ with $\omega^{(i)} \in \Gamma(L,\delta)$ for $i=1,2,\ldots d$ 
\begin{align*}
& \hphantom{\ge}\ \sum_{n \in \N} \prod_{k=1}^{d} P_{\omega^{(k)}}(X_{2n}=0) \vphantom{\sum_{\lceil\exp(3\delta L)\rceil \le n \le \lfloor\exp((1-2\delta)L)\rfloor}} \ge  \sum_{\lceil\exp(3\delta L)\rceil \le n \le \lfloor\exp((1-2\delta)L)\rfloor} \prod_{k=1}^{d} P_{\omega^{(k)}}(X_{2n}=0) \\
&\ge  \Big(\exp\big((1-2\delta)L\big) - \exp(3\delta L) - 1 \Big) \cdot C^d \cdot \exp(-3 \delta d L) \vphantom{\sum^d_{\lceil\exp(3\delta L)\rceil \le n \le \lfloor\exp((1-2\delta)L)\rfloor}}\\
&= C^d \cdot \exp(-3 \delta d L) \cdot \exp(3 \delta d L)   \vphantom{\sum^d_{\lceil\exp(3\delta L)\rceil \le n \le \lfloor\exp((1-2\delta)L)\rfloor}}\allowdisplaybreaks[0]\\
& \qquad \cdot \Big(\exp\big((1-2\delta-3 \delta d)L\big)  - \exp\big((3\delta -3\delta d) L\big) - \exp(-3 \delta d L) \Big) \vphantom{\sum^d_{\lceil\exp(3\delta L)\rceil \le n \le \lfloor\exp((1-2\delta)L)\rfloor}}\\
& \xrightarrow{L\to\infty} \infty \vphantom{\sum^d_{\lceil\exp(3\delta L)\rceil \le n \le \lfloor\exp((1-2\delta)L)\rfloor}}
\end{align*}
for 
\[
0 < \delta < \frac{1}{2+3d}\ .
\]
Since \eqref{I-eq16.2} holds for arbitrarily small $\delta$, we can conclude that \eqref{I-eqcor2} holds for $\p^{\otimes d}$-a.e.\ environment $(\omega^{(1)}, \omega^{(2)}, \ldots, \omega^{(d)})$.
\end{pfof}
\section{Examples for Recurrent Random Walks in Random Environments in different dimensions} \label{I-sec1.6}

\begin{rmk}\label{I-cor1}
Consider a RWRE $(X_n)_{n \in \N_0}$ for which the environment $\omega$ fulfils the assumptions $\eqref{I-ass1}$, $\eqref{I-ass2}$, and $\eqref{I-ass3}$. By an application of Theorem~\ref{I-Rthm1} for $\alpha =0$, we get
\[
\sum_{n \in \N} P_{\omega}(X_{2n}=0) = \infty
\]
for $\p$-a.e.\ environment $\omega$. From this, we can conclude that the random walk is recurrent for  $\p$-a.e.\ environment $\omega$.
\end{rmk}

\begin{cor}[$d$ particles in the same random environment]\label{I-cor5.3a} 
Let us first choose a random environment $\omega=(\omega_x)_{x \in \Z}$ which fulfils the  assumptions $\eqref{I-ass1}$, $\eqref{I-ass2}$, and $\eqref{I-ass3}$. For fixed $\omega$, we can now consider $d$ independent random walks $(X^{(i)}_n)_{n \in \N_0}$ for $i=1,2,\ldots,d$ where every random walk $(X^{(i)}_n)_{n \in \N_0}$ is a usual RWRE in the environment $\omega$ in the sense of \eqref{I-RWRE}. Then, for arbitrary $d$, the $d$-dimensional process
\[
\big(X^{(1)}_n,X^{(2)}_n, \ldots, X^{(d)}_n\big)_{n \in \N_0}
\]
is recurrent for $\p$-a.e.\ environment $\omega$.
\end{cor}

\begin{pfof}{Corollary \ref{I-cor5.3a}}
First of all, we notice that for fixed $\omega$ 
\[
\big(X^{(1)}_n,X^{(2)}_n, \ldots, X^{(d)}_n\big)_{n \in \N_0}
\]
is a Markov chain. For the expected amounts of returns to $0$, we get by applying Theorem \ref{I-Rthm2} with $\alpha=d$
\begin{align*}
\sum_{n\in\N} P_{\omega} \left(\big(X^{(1)}_{2n},X^{(2)}_{2n}, \ldots, X^{(d)}_{2n}\big)= \big(0,0,\ldots,0\big) \right) = \sum_{n\in\N} \left(P_{\omega} (X^{(1)}_{2n} = 0) \right)^d = \infty
\end{align*}
for $\p$-a.e.\ environment $\omega$. This implies the recurrence.
\renewcommand{\qedsymbol}{$\square$}
\end{pfof}

\begin{cor}[$d$ particles in $d$ i.i.d.\ random environments] \label{I-cor5.3}
For arbitrary $d \in \N$, we choose $d$ i.i.d.\ environments $\omega^{(i)}=(\omega^{(i)}_x)_{x \in \Z}$ which all fulfil the assumptions $\eqref{I-ass1}$, $\eqref{I-ass2}$, and $\eqref{I-ass3}$ for $i=1,2,\ldots,d$. For fixed $\vec{\omega}:=(\omega^{(1)}, \omega^{(2)}, \ldots, \omega^{(d)})$, we consider $d$ independent RWRE $(X^{(i)}_n)_{n \in \N_0}$, where $(X^{(i)}_n)_{n \in \N_0}$ is a usual RWRE in the environment $\omega^{(i)}$ in the sense of \eqref{I-RWRE}. In this case, the  $d$-dimensional process
\[
\big(X^{(1)}_n,X^{(2)}_n, \ldots, X^{(d)}_n\big)_{n \in \N_0}
\]
is recurrent for $\p^{\otimes d}$-a.e.\ environment $\vec{\omega}$.
\end{cor}
\begin{pfof}{Corollary \ref{I-cor5.3}}
Due to the independence of the processes and the environments in every component, we get
\begin{align*}
\sum_{n\in\N} P_{\vec{\omega}} \left(\big(X^{(1)}_{2n},X^{(2)}_{2n}, \ldots, X^{(d)}_{2n}\big)= \big(0,0,\ldots,0\big) \right) = \sum_{n\in\N} \prod_{i=1}^{d} P_{\omega^{(i)}} (X^{(i)}_{2n}=0) = \infty
\end{align*}
due to Theorem \ref{I-Rthm3} for $\p^{\otimes d}$-a.e.\ environment $\vec{\omega}$. 
\renewcommand{\qedsymbol}{$\square$}
\end{pfof}

\begin{rmk}
An alternative proof of Corollary \ref{I-cor5.3} can be found in \cite{Zei} after Lemma A.2. The proof there uses the Nash-Williams inequality in the context of electrical networks.
\end{rmk}

\begin{rmk}
Corollary \ref{I-cor5.3a} and \ref{I-cor5.3} show that the recurrence of a RWRE is indeed ``stronger'' than the recurrence of the symmetric random walk on $\mathbb{Z}$. Note that $d$ particles performing a one-dimensional symmetric random walk will only meet finitely often for $d \ge 3$.
\end{rmk}

\begin{cor}[Symmetric Random Walk combined with RWRE - Version 1]\label{I-cor6}
We first choose an environment $\omega$ which fulfils the assumptions \eqref{I-ass1}, \eqref{I-ass2}, and \eqref{I-ass3}. For a fixed $\omega$, let $(X_n,Y_n)_{n \in \N_0}$ be a 2-dimensional process where the process $(X_n)_{n\in \N_0}$ and $(Y_n)_{n\in \N_0}$ are independent with respect to $P_{\omega}$, $(X_n)_{n\in \N_0}$ is a RWRE in the sense of \eqref{I-RWRE} and $(Y_n)_{n\in \N_0}$ a symmetric random walk on $\Z$. Then, $(X_n,Y_n)_{n \in \N_0}$ is recurrent for $\p$-a.e.\ environment $\omega$.
\end{cor}

\begin{pfof}{Corollary \ref{I-cor6}}
Due to the independence of the two components, we get
\begin{align*}
& \sum_{n \in \N} P_{\omega}\big((X_{2n},Y_{2n})=(0,0)\big) = \sum_{n \in \N} P_{\omega}\big(X_{2n}=0\big) \cdot P_{\omega}\big(Y_{2n}=0\big) \\
\ge\ & C \cdot \sum_{n \in \N} P_{\omega}\big((X_{2n}=0\big) \cdot n^{-\frac12} = \infty.
\end{align*}
Here, we used the lower bound
\begin{align} \label{I-cor1.6.6.1}
P_{\omega}\big(Y_{2n}=0\big) \ge C \cdot n^{-\frac12}
\end{align}
for the return probabilities of the symmetric random walk on $\Z$ with some constant $C>0$ (cf.\ Section 2.18.4 in \cite{Gut}) and Theorem \ref{I-Rthm1} with $\alpha =\tfrac12$ for the last two steps. Again, we can conclude the recurrence of the process $(X_n,Y_n)_{n \in \N_0}$ for $\p$-a.e.\ environment $\omega$.
\renewcommand{\qedsymbol}{$\square$}
\end{pfof}

\begin{cor}[Symmetric Random Walk combined with RWRE - Version 2] \label{I-cor1.5.7}
We first choose an environment $\omega$ which fulfils the assumptions \eqref{I-ass1}, \eqref{I-ass2}, and \eqref{I-ass3} and some $0 < \delta <1$. For a fixed environment $\omega$, let $(X_n,Y_n)_{n \in \N_0}$ be a Markov chain with values in $\Z^2$ which is determined by
\begin{align*}
&P_{\omega}\big((X_0,Y_0)=(0,0)\big) = 1 \vphantom{\frac{1 - \delta}{2}},\\
&P_{\omega}\big((X_{n+1},Y_{n+1})=(x+1,y)\big|(X_{n},Y_{n})=(x,y)\big)=\delta \cdot \omega_x \vphantom{\frac{1 - \delta}{2}},\\
&P_{\omega}\big((X_{n+1},Y_{n+1})=(x-1,y)\big|(X_{n},Y_{n})=(x,y)\big)=\delta \cdot (1-\omega_x) \vphantom{\frac{1 - \delta}{2}},\\
&P_{\omega}\big((X_{n+1},Y_{n+1})=(x,y\pm1)\big|(X_{n},Y_{n})=(x,y)\big)=\frac{1 - \delta}{2}\ .
\end{align*}
Again, $(X_n,Y_n)_{n \in \N_0}$ is recurrent for $\p$-a.e.\ environment $\omega$.
\end{cor}

\begin{figure}[h] 
\begin{minipage}[t]{0.35\textwidth}
\begin{rmk}
In the situation of Corollary \ref{I-cor1.5.7}, we first choose the first (or second) component for the next step with probability $\delta$ (or $1-\delta$). If we choose the first component, then we change the first component by $\pm1$ as in the setting of a RWRE, otherwise we change the second component by $\pm1$ with probability $\tfrac12$ as in the case of a symmetric random walk on $\Z$. 
\end{rmk}
\end{minipage}
\hfill
\begin{minipage}[t]{0.5\textwidth}
\vspace{0pt}
\centering
\includegraphics[viewport=320 450 200 800, scale =0.6]{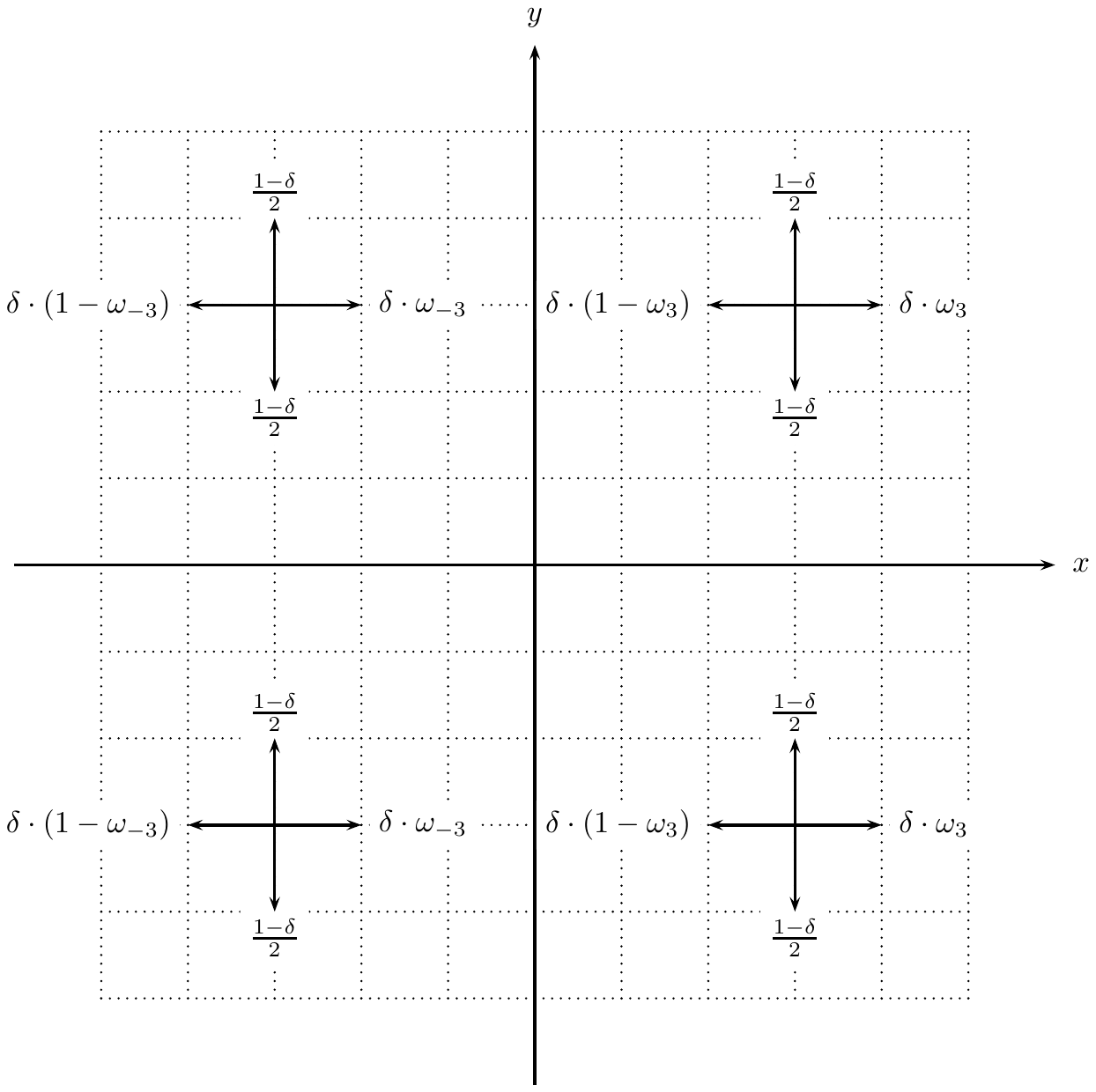}
  \caption{Transition probabilities for the considered process in Corollary \ref{I-cor1.5.7}}
\end{minipage}
\end{figure}

 
\begin{pfof}{Corollary \ref{I-cor1.5.7}}
For the proof, it is enough to look at the process $(X_{n},Y_{n})_{n \in \N_0}$ whenever it has moved in the first component. For this, we define inductively
\begin{align*}
& \tau_0:=0 & \text{and}\\
& \tau_k:= \inf\left\{n > \tau_{k-1}:\ X_n \neq X_{\tau_{k-1}} \right\} & \text{for } k \ge 1.
\end{align*}

\vspace{0pt}
Additionally, we define
\begin{align*}
& \widetilde{X}_n:= X_{\tau_n} \qquad \hphantom{\widehat{Y}_n:= Y_{\tau_n}} \text{for } n \in \N_0,\\
& \widetilde{Y}_n:= Y_{\tau_n} \qquad \hphantom{\widehat{X}_n:= X_{\tau_n}} \text{for } n \in \N_0. 
\end{align*}
Note that $(\widetilde{X}_n)_{n \in \N_0}$ is a usual RWRE on $\Z$ for which the environment $\omega$ fulfils our assumptions \eqref{I-ass1}, \eqref{I-ass2}, and \eqref{I-ass3}. Further, we have
\begin{align} \label{I-cor1.6.7.3}
\widetilde{Y}_n \stackrel{d}{=} S(\tau_n - n),
\end{align}
where $(S(n))_{n \in \N_0}$ denotes a symmetric random walk on $\Z$ which is independent of $(\widetilde{X}_n)_{n \in \N_0}$, $(\tau_n)_{n \in \N_0}$, and the environment $\omega$. Note here that we can decompose $\tau_n$ into the increments
\begin{align}\label{I-cor1.6.7.1}
\tau_n =  \sum_{i=1}^{n} (\tau_i - \tau_{i-1}),
\end{align}
where $(\tau_i - \tau_{i-1})_{i \in \N}$ is a sequence of i.i.d.\ random variables with a geometric distribution with parameter $\delta$ and expectation $\tfrac{1}{\delta}$.

Let us fix an arbitrary $\gamma > 0$. Due to \eqref{I-cor1.6.7.1}, an application of Cramer's theorem implies that we have
\begin{align*}
& P_{\omega} \left(\tau_n > \left(\tfrac{1}{\delta} + \gamma \right) \cdot n \right) \le \exp(-n \cdot I) 
\end{align*}
for some constant $I$=$I(\gamma) > 0$. Therefore, the Borel-Cantelli lemma implies that we have
\begin{align*} 
P_{\omega} \left(\liminf_{n\to\infty} \left\{ n \le \tau_n \le \left(\tfrac{1}{\delta} + \gamma \right) \cdot n \right\} \right) = 1
\end{align*}
for every environment $\omega$. Notice here that we have $\tau_n \ge n$ by definition. Due to the continuity of $P_{\omega}$, we can therefore conclude that
\begin{align} \label{I-cor1.6.7.2a}
\lim_{n \to \infty} P_{\omega} \Big( n \le \tau_n \le \left(\tfrac{1}{\delta} + \gamma \right) \cdot n  \Big) = 1.
\end{align}
Since we are interested in the returns of the random walk to 0, we have to distinguish between the cases in which $\tau_n$ is even or odd. Only for even values of $\tau_n$ our random walk $(\widetilde{X}_{2n}, \widetilde{Y}_{2n}) = (X_{\tau_{2n}}, Y_{\tau_{2n}})$ can reach the point $(0,0)$. For this, we note that $\tau_n$ has a negative binomial distribution with parameters $n$ and $\delta$ and therefore has the following properties:
\begin{align}
& P_{\omega}  (\tau_n = k) \le P_{\omega}  (\tau_n = k + 1) \quad \text{for } n \le k \le \frac{n -1}{\delta} \nonumber,\\
& P_{\omega}  (\tau_n = k) \ge P_{\omega}  (\tau_n = k + 1) \quad \text{for }  k \ge \max \left\{\frac{n - 1}{\delta}, n \right\} \nonumber,\\
& \max_{k \ge n} P_{\omega}  (\tau_n = k) \xrightarrow{n \to \infty} 0. \label{I-cor1.6.7.2b}
\end{align}
Thus, a combination of \eqref{I-cor1.6.7.2a} and \eqref{I-cor1.6.7.2b} implies that in the limit, for $n \to \infty$, the probability for the even and odd part is the same, i.e.
\begin{align*} 
\lim_{n \to \infty} P_{\omega} \Big( n \le \tau_n \le \left(\tfrac{1}{\delta} + \gamma \right) \cdot n  ,\ \tau_n \in 2 \N_0 \Big) = \frac12.
\end{align*}
Since due to our choice $\gamma > 0$ we have
\begin{align}
P_{\omega} \Big(n  \le \tau_n \le \left(\tfrac{1}{\delta} + \gamma \right) \cdot n,\ \tau_n \in 2 \N_0 \Big) > 0 \nonumber
\intertext{for all $n \in \N$, a combination of the last two in-/equalities implies that there exists some constant $C_2 > 0$ such that}
P_{\omega} \Big(n  \le \tau_n \le \left(\tfrac{1}{\delta} + \gamma \right) \cdot n,\ \tau_n \in 2 \N_0  \Big) \ge C_2 > 0\label{I-cor1.6.7.4}
\end{align}
for all $n \in \N$ and for every environment $\omega$. Using the independence of $(\widetilde{X}_{2n})_{n \in \N_0}$ and $(\widetilde{Y}_{2n})_{n \in \N_0}$, we therefore get the following lower bound:
\begin{align*}
& \sum_{n \in \N} P_{\omega}\big((\widetilde{X}_{2n},\widetilde{Y}_{2n})=(0,0)\big) = \sum_{n \in \N} P_{\omega}(\widetilde{X}_{2n}=0) \cdot P_{\omega}(\widetilde{Y}_{2n}=0) \vphantom{\sum_{\substack{i = 2n \\ i \in 2 \N_0}}^{\lfloor \left(\tfrac{1}{\delta}+ \gamma\right) \cdot2 n \rfloor}} \\
\ge\ & \sum_{n \in \N} P_{\omega}(\widetilde{X}_{2n}=0) \cdot \sum_{\substack{i = 2n \\ i \in 2 \N_0}}^{\left\lfloor \left(\tfrac{1}{\delta} + \gamma \right) \cdot2 n \right\rfloor} P_{\omega}(\widetilde{Y}_{2n}=0, \tau_{2n} = i )\\
\ge\ & \sum_{n \in \N} P_{\omega}(\widetilde{X}_{2n}=0) \cdot \sum_{\substack{i = 2n \\ i \in 2 \N_0}}^{\left\lfloor \left(\tfrac{1}{\delta} + \gamma\right) \cdot 2n \right\rfloor} P_{\omega}\big(S(i - 2n) =0 \big) \cdot P_{\omega}(\tau_{2n} = i )\\
\ge \ & \sum_{n \in \N} P_{\omega}(\widetilde{X}_{2n}=0) \cdot \Bigg(P_{\omega}(\tau_{2n} = 2n ) + \sum_{\substack{i = 2n + 2 \\ i \in 2 \N_0}}^{\left\lfloor \left(\tfrac{1}{\delta} + \gamma\right) \cdot 2n \right\rfloor} C \cdot (i - 2n)^{-\frac12} \cdot P_{\omega}(\tau_{2n} = i) \Bigg)
\intertext{Here, we used \eqref{I-cor1.6.7.3} in the third line and the usual lower bound for the return probabilities of the symmetric random walk on $\Z$ (cf.\ \eqref{I-cor1.6.6.1}), i.e. 
\[
P_{\omega}\big(S(i - 2n) =0 \big) \ge C \cdot (i-2n)^{-\frac12}
\]
for $i \in 2\N$, $i \ge 2n+2$ and with some constant $C>0$, in the fourth line. From this, we get}
& \sum_{n \in \N} P_{\omega}\big((\widetilde{X}_{2n},\widetilde{Y}_{2n}=(0,0)\big) \vphantom{\sum_{\substack{i = 2n \\ i \in 2 \Z}}^{\lfloor \left(\tfrac{1}{\delta}+ \gamma\right) \cdot2 n \rfloor}}\\
\ge \ & C \cdot  \sum_{n \in \N} P_{\omega}(\widetilde{X}_{2n}=0) \cdot \big(2 \cdot\left( \tfrac{1}{\delta} + \gamma - 1\right)\big)^{-\frac12} \cdot  n^{-\frac12} \cdot  \sum_{\substack{i = 2n \\ i \in 2 \N_0}}^{\left\lfloor \left(\tfrac{1}{\delta} + \gamma\right) \cdot 2n \right\rfloor} P_{\omega}(\tau_{2n} = i ) \\
\ge\ & C \cdot C_2 \cdot \big(2 \cdot\left( \tfrac{1}{\delta} + \gamma - 1 \right)\big)^{-\frac12} \cdot  \sum_{n \in \N} P_{\omega}(\widetilde{X}_{2n}=0) \cdot n^{-\frac12} = \infty \vphantom{\sum_{\substack{i = 2n \\ i \in 2 \Z}}^{\lfloor \left(\tfrac{1}{\delta}+ \gamma\right) \cdot2 n \rfloor}}
\end{align*}
for $\p$-a.e.\ environment $\omega$. Here, we additionally made use of \eqref{I-cor1.6.7.4} and Theorem \ref{I-Rthm1} (applied for $\alpha =\tfrac12$) in the last line. This implies that the process 
\[
(\widetilde{X}_{n},\widetilde{Y}_{n})_{n \in \N_0}
\]
is recurrent for $\p$-a.e.\ environment $\omega$. Finally, this obviously implies that our process
\[
(X_{n},Y_{n})_{n \in \N_0}
\]
is also recurrent for $\p$-a.e.\ environment $\omega$ since we can embed the paths of the process $(\hspace{-1pt}\widetilde{X}_{n},\hspace{-2pt}\widetilde{Y}_{n}\hspace{-1pt})_{n \in \N_0}$ into the paths of the process $(X_{n},Y_{n})_{n \in \N_0}$.
\renewcommand{\qedsymbol}{$\square$}
\end{pfof}

\newpage
 
\bibliographystyle{alpha}

\bigskip\bigskip\bigskip
\noindent
Michael Kochler\\
Center for Mathematical Sciences\\
Technische Universit\"at M\"unchen\\
Boltzmannstra\ss e 3 \\
D-85748 Garching\\
Germany\\
\noindent e-mail: \texttt{michael.kochler@tum.de}

\end{document}